\documentclass{amsproc}
\usepackage{amsfonts}
\usepackage{amsmath}
\usepackage{enumerate}
\usepackage{graphicx}

\newtheorem{theorem}{Theorem}

\newtheorem{definition}[theorem]{Definition}

\newcommand{\Pp}{\mathbf{P}}
\newcommand{\R}{\mathbb{R}}
\newcommand{\N}{\mathbb{N}}

\newcommand{\ONE}{{\bf 1}}
\newcommand{\eps}{\varepsilon}

\newcommand{\x}{X_\eps}

\begin{document}
\title{Scaling Limit for the Diffusion Exit Problem}

\author{Sergio A.Almada}

\address{Department of Statistics and Operations Research, University of North Carolina, Chapel Hill}
\curraddr{UNC Department of Statistics \& Operations Research, 318 Hanes Hall, CB\# 3260, UNC Chapel Hill, NC 27599-3260}

\email{salmada3@email.unc.edu}


\begin{abstract} In this review, an outline of the so called Freidlin-Wentzell theory and its recent extensions is given. Broadly, this theory studies the exponential rate at which the probabilities of rare events related to random perturbation of ODE decays. The typical situation is when an
ODE has several stable equilibria, in which case, the theory predicts the most likely paths in which the randomly perturbed system goes from one equilibria to another. In recent developments I will outline how recent approaches allows to distinguish between paths that are otherwise exponentially equivalent. An overview of applications of this theory is briefly covered.
\end{abstract}
\maketitle
\section*{Introduction}
In this survey we study the so called exit problem~\cite[Section 4.3]{Freidlin--Wentzell-book} for small noise diffusion. This model belongs to the more general area of random perturbations of dynamical systems, which has been a very active area of research over the last 30 years~\cite{GentzBook},~\cite{Freidlin--Wentzell-book},~\cite{Vares}. The small noise diffusion framework has attracted  the interest of both the  pure and applied mathematics communities. From the mathematical standpoint  it is interesting partly
because this area has strong interactions with other important branches of mathematics such as probability theory, dynamical systems, or PDE. As regards applied mathematics,  the set of problems relating to small noise diffusion has found applications in climate modeling~\cite{Clima},~\cite{GentzClimate}, electrical engineering~\cite{Zei2},~\cite{Zei1},~\cite{Zei3},  finance~\cite{FinancePaper},~\cite{FinanceBook}, neural dynamics~\cite{rabinovich:188103},~\cite{Biology} among others~\cite{DemboZeitouni}. 

The setting of the problem is the following. Given a smooth vector field $b: \R^d \to \R^d$ consider the It\^o equation driven by the $d$-dimensional standard Wiener process $W$:
\begin{align}
\label {eqn: Ito_additive} d\x (t) &= b( \x (t) ) dt + \eps dW(t), \\
\x (0 ) &= x_0. \notag
\end{align}
Assume that the vector field $b$ is such that we can ensure that equation~\eqref{eqn: Ito_additive} has a unique strong solution (see~\cite{Karatzas--Shreve} or~\cite{Protter} for all stochastic analysis references). Given an initial condition $x_0 \in \R^d$ (or a set of initial conditions), the goal is to characterize some asymptotic properties of $\x$ as $\eps \to 0$. 

The focus of this survey is on the exit from a domain problem or exit problem for short. Consider a domain (open, bounded and connected) $D \subset \R^d$ with piecewise smooth boundary (at least $C^2$). The exit problem is the study of the time
\[
\tau_\eps^D (x)= \inf \{ t>0: \x(t) \in \partial D \},
\]
at which $\x$ exits the domain $D$, and the exit distribution $\Pp_{x_0} \{ \x ( \tau_\eps^D ) \in \cdot \} $. As is expected the asymptotic distribution depends on the dynamic properties inherited by the flow generated by $b$. Some of this properties have been leveraged to several applications that we will briefly mention. 

	The standard mindset in tackling this problem from the probabilistic point of view is to think of the SDE that defines $\x$ as a (random) singular perturbation of the system $\dot{x}=b(x)$. Freidlin and Wentzell~\cite{Freidlin--Wentzell-book},~\cite{Vares} were the ones who put together a general theory in this direction, mostly based on the Large Deviation principle for $\x$ (see Section~\ref{sec: app-LDP} for a brief review on Large Deviations). The theory came to light with a series of papers beginning with~\cite{FWPaper1} and~\cite{FWPaper2} until the Russian edition of the book~\cite{Freidlin--Wentzell-book} appeared.  See~\cite{KoralovSPA} and ~\cite{KoralovPTRF} for a modern version of the theory, and~\cite{SPDE1},~\cite{Cerra1} and references therein for a stochastic partial differential equations version of the theory. In Section~\ref{sec: intro_FW} we give a brief review of the Freidlin-Wentzell theory.
	
In contrast with the Freidlin-Wentzell theory, that mostly relies on the large deviation principle, a modern trend relying on a path-wise approach has emerged in the last years. As a consequence, more detailed phenomena can be captured. That is the case, for example, in~\cite{nhn} in which a heteroclinic network is considered or in~\cite{GentzP1} in which a bifurcation problem is studied. The monograph~\cite{GentzBook} contains several examples in this direction together with applications.

In this survey, we also present recent applications and developments to Monte Carlo algorithms~\cite{MonteCarlo}. The connection between the exit problem and these algorithms is established using the statistical physics framework~\cite{GentzQuimica},~\cite{FreeEnergy}. In particular we focus on the Simulated Annealing algorithm, ~\cite{TermodinamicAnnealing},~\cite{SimulatedAnnealing}, which is one of the most used algorithms in practice. From a mathematical point of view, this algorithm was studied in a series of papers~\cite{DiffusionRn},~\cite{Gelfand},~\cite{Geman},~\cite{Kushner}, and~\cite{AnnealingRates}.

Section~\ref{sec: intro_FW} contains the basic setup and results. This section assumes basic intuitive understanding of Large Deviations Theory, if the reader lacks such understanding, Section~\ref{sec: app-LDP} is meant to be self contained and to fill this gap. Section~\ref{sec: Saddle} contains the case of a saddle point. Section~\ref{sec: Applications} contains an intuitive explanation of an application (Simulated Annealing) that has been influenced by the small noise problem. Section~\ref{sec: FutureDirections} briefly presents a view on future directions. 
 
\section{Background and Motivation}\label{sec: intro_FW}

Let $\x$ be the strong solution to the SDE~\eqref{eqn: Ito_additive}. The equation for $\x$ suggests that the process should, for small $\eps$, behave like the flow generated by $b$:
\begin{align}
\label {eqn: det_flow} \frac{d}{dt}S^t x = b(S^tx ), \quad
S^tx = x.
\end{align}
Indeed, through a standard martingale argument, it is easy to see that for any $\delta>0$ there are constants $C_1=C_1({T,\delta})$ and $C_2({T,\delta})$ such that
\begin{equation} \label{eqn: estimate_basic}
\sup_{ x \in \R^d} \Pp_x \left \{  \sup_{ t \leq T} | \x (t) - S^t x| > \delta \right \} \leq C_1 e^{ - C_2\eps^{-2} }.
\end{equation}
Note that this estimate applies only on a compact time interval $[0,T]$. In principle $\tau_\eps$ can grow  to infinity asymptotically for small $\eps$, which limits the usability of~\eqref{eqn: estimate_basic}. The first building block to overcome this difficulty is to establish a Large Deviation Principle (LDP) for $\x$. In this case, finding the LDP for $\x$ consists on finding the optimal constant $C_2$ in~\eqref{eqn: estimate_basic}, see~\cite{DemboZeitouni},~\cite{DenHollanderLDP} or Section~\ref{sec: app-LDP}. The LDP for $\x$ is the following:
\begin{theorem} [Freidlin-Wentzell~\cite{Freidlin--Wentzell-book} ] \label{thm: LDP}
Let $H_{0,T}^1$ be the space of all absolutely continuous functions from $[0,T]$ to $\R^d$ with square integrable derivatives. Define the functional $I_T^x$ by
\begin{equation*}
I_T^x (\varphi)= \frac{1}{2} \int_0^T  \| \stackrel{\cdot}{\varphi} (s)-b(\varphi(s) ) \| ds,
\end{equation*}
if $\varphi \in H_{0,T}^1$ and $\varphi(0)=x$, and $\infty$ otherwise. 

Then for each $x \in \R^d$ and $T>0$ the family $(\Pp_x^\eps)_{\eps>0}$ satisfies a Large Deviation Principle on $C([0,T];\R^d)$ equipped with uniform norm at rate $\eps^2$ with good rate function $I_T^x$.
\end{theorem}

Informally, Theorem~\ref{thm: LDP} says that if $A \subset C([0,T];\R^d)$ then
\begin{equation}  \label{eqn: intro-LDP_I}
-\eps^2 \log \Pp_x \left\{   \x \in A \right\} \asymp \inf_{\varphi \in A} I_T^x (\varphi), \text{ as } \eps \to 0,
\end{equation}
which implies that the optimal constant $C_2$ in~\eqref{eqn: estimate_basic} is given by 
$C_2 = T \delta/2.$ Further,~\eqref{eqn: intro-LDP_I} suggests that $I_T^x$ can be viewed  as a measure on how costly (in terms of probability) is for the system $\x$ not to follow the deterministic trajectory $S$. This interpretation is essential when solving problems that require non-compact time frames.

Regarding $I_T^x$ as a cost function, it make sense to introduce $V: D \times  \partial D \to [0, \infty]$ given by
\begin{equation} \label{eqn: quasipotential}
V(x,y) = \inf_{T>0} \left \{ I_T^x ( \varphi) : \varphi (T)=y, \varphi( [0,T] ) \subset D \cup \partial D \right \},
\end{equation}
which represents the cost that the process $\x$ would incur to go from $x$ to $y \in \partial D$. This function is known as the quasi-potential, and it plays an important role on the exit problem:

\begin{theorem} [Freidlin-Wentzell~\cite{Freidlin--Wentzell-book}]
Suppose $b$ is smooth, $\x (0)=x_0$ and let \[z=\inf_{y\in \partial D } V(x_0,y).\] Then, there is a constant $\bar V$ such that for every $\delta \in [0, \bar V)$ and for every closed set $N \subset \partial D$ that satisfies $\inf_{y\in N } V(x_0,y) > z$,
\[
\lim_{\eps \to 0} e^{ ( \bar V - \delta )/ \eps^2 } \Pp_{x_0} \{ \x (\tau_\eps) \in N \}=0.
\]
\end{theorem}

The immediate observation resulting from this theorem is that the location of the exit distribution is concentrated on the set of minimizers $V_*$ of the quasipotential.  At the same time, it implies that for a set $N \subset \partial D$ that contains $V_*$, 
\[
\Pp_{x_0} \{ \x (\tau_\eps)  \not \in N \} = p_\eps e^{- \bar V/\eps^2 },
\]
where $\eps^2 \log p_\eps \to 0$. As such, in principle, this theorem in itself does not provide any information on the exit distribution when the exit is restricted to $V_*$. This kind of information is usually obtained via ad-hoc analysis to the particular case under consideration. In the next section we present one of these cases.

\section{Saddle Case, an Assymmetric Example} \label{sec: Saddle}
In this section we consider the small noise scape from a saddle problem that, up to our knowledge, was first studied in~\cite{Kifer}. Our objective is to ilustrate the results obtained in cases in which the quasi-potential approach provides incomplete information. For simplicity we will focus on the case in which $d=2$. 

The saddle case can be described as the case in which $0 \in D$ is the only critical point of $b$ in the closure of $D$; that is, $0 \in D$ is the only $x\in \bar{D}$ such that $b(x)=0$. Further, suppose that the vector field $b$ is such that its Jacobian at $0$, $A=Db(0)$ has at one positive eigenvalue and one negative eigenvalue. The case of interest for this problem is when the initial condition for the diffusion $X_\eps(0)$ lies in the invariant stable manifold $$\mathcal{M}^s = \left \{  x: S^t x \to 0, \text{ as } t \to \infty \right \}.$$ In this case, it can be shown that the quasi-potential has two minimizers $\{ q_-, q_+ \}$ that correspond to the two intersection points of the unstable invariant manifold 
\[
 \mathcal{M}^u =  \left \{  x: S^t x \to 0, \text{ as } t \to -\infty \right \}
 \]
 with $\partial D$. As seen on the last section, the exit distrition of $\x$ from $D$ is concentrated on these two points. in this section we will study the distribution precisely.
 
The small noise exit problem from a saddle point was first solved using a PDE approach in~\cite{Kifer}. In that paper, it is shown that the exit time is asymptotically logarithmic in $\eps$ and that the exit location is uniform on the set $\{ q_-, q_+ \}$. Later,~\cite{Day} refined the result of the exit distribution in two dimensions, and further refinements were made in higher dimensions in~\cite{Bakhtin-SPA}.

In~\cite{nhn} a further generalization to the exit location was obtained, and it was shown that if a small perturbation of the initial condition is applied, then an asymmetric exit distribution is obtained. The final result in this two dimensional setting was obtained in~\cite{MioNonl}, which is were the formulation of the following theorem was taken.

\begin{theorem} \label{thm: SmallNoise}
Suppose that $A$ has spectrum $\lambda_+ > 0 > - \lambda_-$ and denote $\partial U \cap \mathcal{M}^u = \{q_-, q_+\}$, and assume that $X_\eps(0) = x_0 +\eps^\alpha v$, with $x_0 \in \mathcal{M}^s \cap U, \alpha > 0$, and $v \in \R^2$. 

Then, there is a family of random vectors $( \phi_\eps )_{\eps >0}$, a family of random variables $( \xi_\eps )_{ \eps > 0 }$, and a number  
\[
\beta = \left\{
	\begin{array}{ll}
		1,  & \mbox{if }  \alpha \lambda_- \geq \lambda_+, \\
		\alpha \frac{\lambda_-}{\lambda_+}, & \mbox{if } \alpha \lambda_- < \lambda_+ \\ 
	\end{array}
\right.
\]
such that $X_\eps( \tau_\eps ) =q_{ { \rm sgn } ( \xi_\eps )} + \eps^\beta \phi_\eps$, and the random vector 
\[
\Pi_\eps = \left( \xi_\eps, \phi_\eps, \tau_\eps + \frac{\alpha}{\lambda_+} \log \eps \right)
\]
converges in distribution as $\eps \to 0$.
\end{theorem}
As seen from this theorem, the choice of $q_-$ or $q_+$ depends on the sign of the random variable $\xi_\eps$. To be more precise, the asymptotic exit law of the process is fully determined from the distribution of ${\rm sgn} ( \xi_\eps )$ which, as we will illustrate in the following, depends on the perturbation $v$. This distribution where studied in detail in~\cite{MioNonl}, and~\cite{nhn}.

To give a survey of how this asymmetry is created, let us study the result stated in Theorem~\ref{thm: SmallNoise} in a particular case. Suppose, $b(x,y)=( \lambda_+ x, - \lambda_- y), D=[-1,1]\times [-1,1], x_0=(0,y_0)$, and $v=(\nu, 0)$. In this case, $q_\pm = ( \pm1,0 )$ and everything else is solvable making the proofs from~\cite{MioNonl}, and~\cite{nhn} easy to follow and apply to this case. Some sample computations in a similar case are performed in the survey~\cite{BakhtinSurvey}. We will summarize the results, and restrict our selfs to the case in which $\alpha \in (0,1]$, since its when the asymmetry exist, still we will give the final result in the case in which $\alpha >1$ for completeness and discussion. 

Let us define the random variable 
\begin{equation} \label{eqn: Xi0}
\xi_0 =\nu + \ONE_{ \left( \alpha = 1 \right) } \int_0^\infty e^{-\lambda_+ s }dW_1(s),
\end{equation}
then, $\xi_\eps \to \xi_0$ in probability as $\eps \to 0$. Further, in the case that $\alpha \lambda_- < \lambda_+$, 
\[
\Pi_\eps \stackrel{\Pp}{\longrightarrow} \left( \xi_0, y_0 \left| \xi_0 \right|^{\lambda_- / \lambda_+ }, -\frac{1}{\lambda_+} \log | \xi_0 | \right), \text{ as } \eps \to 0,
\]  
while for $\alpha \lambda_- \geq \lambda_+$, $\Pi_\eps$ converges in distribution to a random variable $\Pi_0$ that has the same distribution as 
\[
\left( \xi_0, y_0 \ONE_{\left(\alpha \lambda_- = \lambda_+ \right)}\left| \xi_0 \right|^{\lambda_- / \lambda_+ } + \mathcal{N}, -\frac{1}{\lambda_+} \log | \xi_0 | \right),
\]
where $\mathcal{N}$ is a zero mean normally distributed random variable with variance $1/(2\lambda_-)$ independent of $\xi_0$. By analyzing~\eqref{eqn: Xi0} we can conclude that as $\eps \to 0$, the asymptotic exit law, $$\eta_0( \cdot ) = \lim_{\eps \to 0} \Pp_{x_0} \left\{ \x(\tau_\eps ) \in \cdot \right\},$$ is such that
\begin{enumerate}
\item if $\alpha < 1$, then $\eta_0 \left( \{ q_{ {\rm sgn} (\nu) } \} \right) = 1$;
\item if $\alpha = 1$, then $\eta_0 \left( \{ q_{ {\rm sgn} (\nu) } \} \right) > 1/2$, and the mean of $\eta_0$ is $q_{ {\rm sgn} \nu } $; and,
\item if $\alpha > 1$, then $\eta_0$ is uniform on $\left\{ q_-, q_+\right\}$.
\end{enumerate}

\begin{figure}
\centering
 \includegraphics[height = 2.5in]{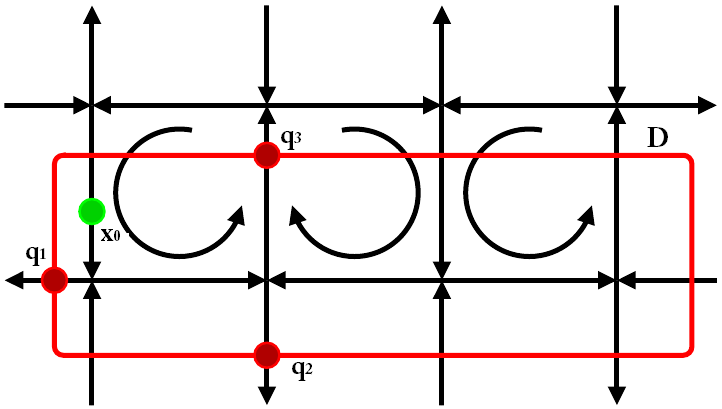}  
 \caption{An example of a noisy heteroclinic network}
 \label{fig: Lattice}
\end{figure} 

The consequences of this asymmetric behavior were applied in~\cite{nhn}, to study the heteroclinic network case. The heteroclinic network case, is the case in which the vector field $b$ has a finite set of critical points $\{ x_0, ..., x_r \}$ such that each point is a saddle, and further the unstable manifold of one critical point is the stable manifold of another critical point forming a network, see Figure~\ref{fig: Lattice}. In this case, if we start in the stable manifold of one of the critical points, then the exit from a neighborhood of this critical point, will be concentrated on the unstable manifold, which is the stable manifold of another critical point. Iterating the procedure again we see that we are now in the situation of Theorem~\ref{thm: SmallNoise} where $\alpha$ is the coefficient $\beta$ from the previous critical point. Hence, an asymmetry can be created. For concreteness, this implies that in figure~\ref{fig: Lattice} the chance of the exit happening on the neighborhood of $q_2$ might be different from the exit happening on the neighborhood of $q_3$ provided that the starting point is $x_0$. This result was unexpected given the history of the problem in question, and the phenomena was first discovered in~\cite{nhn}. We point the reader to~\cite{BakhtinSurvey} for further explanation of this case.

\section{Applications and an informal view of metastability}\label{sec: Applications}

A standard model in statistical physics is to use Gibbs measures to model the state space of our phenomena. Indeed, it is standard to assume that our different states can be characterized by vectors in $\R^d$, and that states can be sampled from the probability distribution
\[
\mu_\eps (dx) = \frac{e ^{- U(x)/ \eps^2 }}{Z_\eps} dx.
\]
In this case, $U(x)$ plays the role of energy of the configuration while $\eps^2$ is proportional to the inverse temperature, and $Z_\eps$ is the free energy that makes $\mu_\eps$ a probability measure: 
\[
Z_\eps = \int_{\R^d} e ^{- U(x)/ \eps^2 }dx .
\]
The connection between this setting, and our standard stochastic differential equation setting is the fact that, under very general conditions, $\mu_\eps$ is the invariant measure of the law of the process $\x$ when $b= -\nabla U$; that is, if the initial condition $\x(0)$ is distributed as $\mu_\eps$, then so is $\x(t)$ for every $t >0$.  In the following $U$ will be a smooth convex function that tends to infinity as $|x| \to \infty$.

This connection between small noise stochastic differential equations and simulation has attracted lots of attention in the applied community. It has a strong connection not only in simulation as explained above (see~\cite[Chapter 9]{MonteCarloBook}  or~\cite{FreeEnergy} for a detailed review) but also with global optimization methods and algorithms. Indeed, the link between this setting and global optimization can be seen from the following theorem which is a summary of the results found in~\cite{LaplaceHwang}: 
\begin{theorem} \label{thm: Hwang}
Suppose that there is a $\eta > 0$ such that $\left \{ x: U(x) \leq \eta \right \}$ is compact. Then, $( \mu_\eps )_{\eps > 0}$ forms a tight family of probability measures, and 
\[
U_* = \left\{ x: x \text{ is a minimizer of } U  \right\}
\]
is not empty. Further, if $\mu_\eps \to \mu$ weakly,  then $\mu$ concentrates on $U_*$.
\end{theorem}

Theorem~\ref{thm: Hwang} makes this setting very attractive for applications for several reasons. For starters from the simulation point of view, allows to simulate discrete objects (represented by a point in $\R^d$) and complex combinatorial relationships among them (encoded in the energy function $U$) by simply discretizing a stochastic differential equation. Moreover, the result of the simulation is likely to correspond to a minimum energy state, which is a common assumption when dealing with physical systems. From the optimization point of view, it in principle allows to find a minimizer of the function $U$ by numerically solving a stochastic differential equation and evolve this solution for a long time. 

As expected this way of thinking has several set backs. In particular, one of the most influential and most studied is the metastability phenomena~\cite{Freidlin--Wentzell-book}. Informally, the process $\x$ is on a metastable state if it spends an exponentially large amount of time in that state. This phenomena is known to happen in cases where $U$ has multiple local minima. 

It is apparent from our discussion that a local minima of $U$, corresponds to a critical point for $b = - \nabla U$. Hence, to study metastable phenomena we need to study the way in which the stochastic process $\x$ passes from a neighborhood of a critical point to a neighborhood of another critical point. Large Deviations Theory in this case provides a well studied answer that builds upon the following well known result.

\begin{theorem} \label{thm: FW}
Suppose $x_*$ is a stable critical point of the driving vector field $b$, and $x_0$ belong to the basin of attraction $B(x_*)$ of $x_*$. That is, $b(x_*) = 0$, and $x_0 \in B(x_*) = \left \{ x: \lim_{t \to \infty } S^tx = x_* \right \}$.

If $D \subset B(x_*)$  strictly, then the first exit time of $D$ is such that 
\[
\lim_{ \eps \to 0 } \Pp_{x_0} \left\{ e^{ (V(x_*) -\delta)/\eps^2 } \leq \tau_\eps^D \leq e^{ (V(x_*) + \delta)/\eps^2 }\right\} = 1,
\]
where $V(x_*) = \inf_{ y \in \partial D } V(x_*, y )$, and $V$ is the quasi potential from Section~\ref{sec: intro_FW}.
\end{theorem}

\begin{figure}
 \centering
 \includegraphics[height = 2.5in]{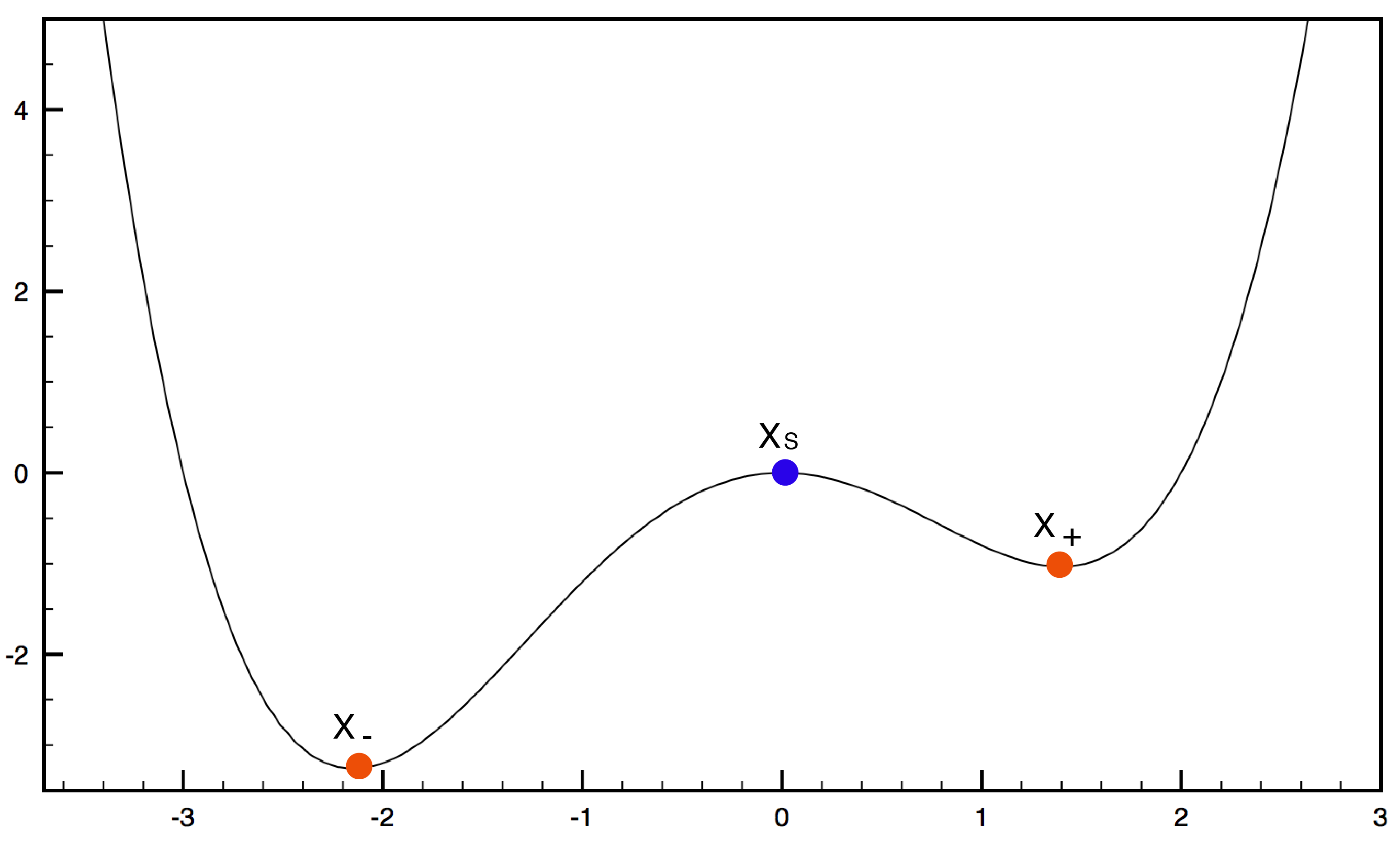}  
 \caption{An example of a double well potential.}
 \label{fig: DoubleWell}
\end{figure}

This theorem allows to study the time that it takes for the process $\x$ to go from one critical point to another by just looking at the quasi potential $V$. This is today a well developed mathematical theory which can be found in the classical reference~\cite{Freidlin--Wentzell-book} or in more modern ones like~\cite{Vares}. For the sake of discussion, suppose we are in a situation where $U$ is a double well potential as depicted in Figure~\ref{fig: DoubleWell}: $b=-\nabla U$ has three critical points, two asymptotically stable $x_-$, and $x_+$, and one saddle point $x_s$ that separates them. In this case, it is well known~\cite{Freidlin--Wentzell-book}[Chapter 4] that $V(x_\pm) = 2\left( U(x_s) - U(x_\pm) \right)$, and Theorem~\ref{thm: FW} implies that in order to escape from the basin of attraction of point $x_\pm$, and fall into the basin of attraction of $x_\mp$ we need to wait a time that is approximately $e^{ 2( U( x_s ) - U( x_\pm ) )/\eps^2 }$. In this case, If we are in a situation in which $x_0$ is on the basin of attraction of $x_-$, and want the process to visit the basin of attraction of $x_+$, we would be stuck for an exponentially long amount of time in the basin of attraction of $x_-$, making the simulation potentially unfeasible. 

In the optimization setting implied by Theorem~\ref{thm: Hwang}, Theorem~\ref{thm: FW} says that is possible to spend an exponentially long amount of time exploring an undesirable minimum: for example, in the double well potential case depicted in Figure~\ref{fig: DoubleWell}, if $x_0 \in B(x_+)$, and $x_-$ is a global minimum of $U$. Several alternatives have been proposed to overcome this issue in practice. Few of those have reached a mathematical level of maturity. Among the ones that have, the most successful seems to be the so called Simulated Annealing~\cite{TermodinamicAnnealing},~\cite{SimulatedAnnealing}.

Simulated Annealing is an applied technique based on the well established Metropolis algorithm~\cite{Metropolis}. It was proposed as an energy relaxation technique in ~\cite{TermodinamicAnnealing}, and~\cite{SimulatedAnnealing}, and ever since has been widely applied in all areas of science and engineering. In our setting, this technique can be vaguely understood based on the observation that for large values of $\eps$ the process $\x$ tends to transition from one basin of attraction to another in a more common way. So the idea is to allow the diffusion to transition easily from critical point (or local minimizer) to critical point by having a relatively large $\eps$ and then reduce the value of $\eps$ once we have reach a candidate critical point.  One way to achieve this effect is by allowing the quadratic variation term to depend on time and to vanish in the long run:
\[
dZ_\eps (t) = - \nabla U ( Z_\eps (t) )dt + \eps(t) dW(t), 
\] 
where $\eps(t) \to 0$, as $t \to \infty$. In a series of papers~\cite{DiffusionRn},~\cite{Gelfand},~\cite{Geman},~\cite{Kushner} have obtained the function form that $\eps(t)$ needs to have in order to establish convergence of $Z_\eps$ to a minimizer of $U$. Their results are summarized in the following theorem, which proof can be found in~\cite{AnnealingRates}.

\begin{theorem}
Suppose $\eps(t)^2 = c/\ln(t)$, and $U$ is a smooth vector field such that $U(z) \to \infty$ as $z \to \infty$, then for every compact set $\Gamma \subset \R^d$ and $c>0$ large enough,
\[
\lim_{ t \to \infty } \sup_{ x_0 \in \Gamma } \eps(t)^2 \ln \Pp_{x_0} \left\{ U(Z_\eps(t) ) \geq \inf_{ z \in \R^d} U(z) +r  \right\} = -2r.
\]
\end{theorem}
This theorem not only establishes the fact that the diffusion $Z_\eps$ will eventually reach a neighborhood of an absolute minima, but it also says how much time it will take to reach it.

\section{Final Remarks and Future Directions} \label{sec: FutureDirections}
As it is clear, algorithms to explore landscapes (energy profiles $U$) have to be motivated by the understanding of the transition process between stable critical points per se. As such, it would be a fundamental tool to develop the necessary theory to understand transitions from one critical point to another similar to the results outlined in Section~\ref{sec: Saddle}. Several advances have been made in this direction. 

A computational theory on how the transitions among basins of attraction of different stable critical points occur is well developed, for example, in~\cite{TransitionPathTheory}. This set of tools heavily rely on the concept of a reactive path: the path that the process $\x$ follows when it scales from the basin of attraction of one stable critical point into the basin of attraction of another stable critical point. As it has been recently pointed out~\cite{JianfengReactivePaths}, a technique to understand these trajectories is via conditioning: conditioned on transition to the basin of a attraction of the stable critical point $A$ how does scape from the basin of attraction of another stable critical point is characterized. The idea of conditioning has had several uses in the past before. In~\cite{MioSPA}, and~\cite{MioMultiScale} the conditioning technique was used to get precise estimates for the exit problem in the simpler case in which no critical points are inside of $D$. More importantly, in~\cite{dayConditional} the idea of conditioning was combined with the quasi-potential framework to study the exit problem from a stable critical point in the case $D$ contains the critical boundary; that is, $\partial D$ is not contained in the basin of attraction of the critical point. On the same flavor, in~\cite{Stein} a PDE approach was used to discover an asymmetric behavior similar to the one described in Section~\ref{sec: Saddle}, but for the characteristic boundary case. It seems that these two results,~\cite{dayConditional} and ~\cite{Stein}, combined with the idea of conditioning could provide the rigorous arguments to study border-line cases in the computational Transition Path Theory~\cite{TransitionPathTheory}. Indeed, this theory relies mostly on estimates derived from Theorem~\ref{thm: FW}, and as such they ignore finer events that are exponentially equivalent as in the case of Section~\ref{sec: Saddle}.

Another approach put forward in the applied community is to modify the idea of Simulated Annealing: instead of letting the quadratic variation term being time dependent, we modify the drift term. The objective of the modification is to alter the potential $U$ so that already visited neighborhoods are easier to scape than unvisited neighborhoods. This type of algorithms are usually known as meta-dynamics, see~\cite{Metadynamics}, and, for a recent survey,~\cite{MetadynamicsSurvey}. In this case, theory is sparse and the formal mathematical properties of these class of algorithms are to be determined. The main difficulty is that the mathematical analysis has to rely on Stochastic Differential Equations with memory. The general large deviation theory for this kind of equations is not well studied up to the authors knowledge. 

\section{Appendix: Large Deviations Principle (LDP) } \label{sec: app-LDP}

Let $\mathcal{X}$ be a Polish metric space with metric function $d:\mathcal{X} \times \mathcal{X} \to [0,\infty)$. By a probability measure on $\mathcal{X}$, we mean a probability measure on the Borel sigma algebra on $\mathcal{X}$. We will give the general definition of large deviation principle for a family of probability measures on $\mathcal{X}$. First, recall the following definition.
\begin{definition} The function $f: \mathcal{X} \to [-\infty,\infty]$ is lower semi-continuous if it satisfies any of the following equivalent properties:
\begin{enumerate}
\item $\liminf_{ n \to \infty } f (x_n) \geq f(x)$ for all sequences $(x_n)_{n\in \N} \subset \mathcal{X}$ and all points $x \in \mathcal{X}$ such that $x_n \to x$ in $\mathcal{X}$.
\item For all $x \in \mathcal{X}$, $\lim_{ \delta \to 0} \inf_{ y \in B_\delta (x) } f(y)=f(x)$, where $B_\delta(x)=\{ y\in \mathcal{X}: d(x,y)<\delta \}$.
\item $f$ has closed level sets, that is, $f^{-1}( [ -\infty,c] )=\{ x \in \mathcal{X}: f(x) \leq c  \}$ is closed for all $c \in \R$.
\end{enumerate}
\end{definition}

Here are the key definitions of large deviation theory:
\begin{definition}
The function $I:\mathcal{X}\to[0,\infty]$ is called a rate function if 
\begin{enumerate}
\item $I \not \equiv \infty$,
\item $I$ is lower semi-continuous,
\item $I$ has compact level sets.
\end{enumerate}
\end{definition}

\begin{definition} \label{def: app-LDP}
A family of probability measures $(\Pp_\eps)_{\eps>0}$ on $\mathcal{X}$ is said to satisfy , as $\eps \to 0$,the large deviation principle (LDP) with rate $\alpha_\eps \to 0$ and rate function $I$ if
\begin{enumerate}
\item $I$ is a rate function,
\item $\limsup_{\eps \to 0 } \alpha_\eps \log \Pp_\eps ( C ) \leq -I ( C )$, for every $C \subset \mathcal{X}$ closed,
\item $\liminf_{\eps \to 0 } \alpha_\eps \log \Pp_\eps ( O ) \geq -I ( O )$, for every $O \subset \mathcal{X}$ open.
\end{enumerate}
Here the bounds are in terms of the set function defined by
\[
I(S)=\inf_{ s \in S } I(x), \quad S \subset \mathcal{X}.
\]
\end{definition}

The goal of large deviation theory is to build up an arsenal of theorems based on these two definitions. We will not describe most of this theorems, since they are out of the scope for the present text. 

\bibliographystyle{plain}
\bibliography{Survey}
\end{document}